\documentclass[11pt]{amsart}
\usepackage{amsbsy,amssymb,amsmath,amsthm,amscd,amsfonts,latexsym,amstext,delarray,
amsmath,graphicx} 
\usepackage[margin=1in]{geometry}
\usepackage{color}
\input xypic

\newtheorem{thm}{Theorem}[section]
\newtheorem{prop}[thm]{Proposition}
\newtheorem{cor}[thm]{Corollary}
\newtheorem{lem}[thm]{Lemma}

\newtheorem{cond}[thm]{Condition}
\newtheorem{defn}[thm]{Definition}
\newtheorem{rem}[thm]{Remark}
\newtheorem{ex}[thm]{Example}

\numberwithin{equation}{section}

\def\F{{\mathbb F}}
\def\Q{{\mathbb Q}}
\def\Z{{\mathbb Z}}
\def\N{{\mathbb N}}

\def\bG{{\mathbb G}}
\def\P{{\mathbb P}}
\def\A{{\mathbb A}}
\def\GL{{\rm GL}}
\def\bL{{\mathbb L}}

\def\cC{{\mathcal C}}

\def\cI{{\mathcal I}}
\def\cJ{{\mathcal J}}

\def\cM{{\mathcal M}}

\def\cO{{\mathcal O}}
\def\cP{{\mathcal P}}

\def\cT{{\mathcal T}}

\def\cV{{\mathcal V}}

\def\cX{{\mathcal X}}
\def\cY{{\mathcal Y}}

\def\bA{{\mathbb A}}

\def\bC{{\mathbb C}}

\def\bG{{\mathbb G}}

\def\bL{{\mathbb L}}

\def\bT{{\mathbb T}}

\title{$\F_\zeta$-geometry, Tate motives, and the Habiro ring}
\author{Catharine Wing Kwan Lo and Matilde Marcolli}
\address{Mathematics Department, Caltech, 1200 E. California Blvd. Pasadena, CA 91125, USA}
\email{cwlo@caltech.edu}
\email{matilde@caltech.edu}
\address{Physics Department, Caltech, 1200 E. California Blvd. Pasadena, CA 91125, USA}
\email{cwlo@caltech.edu}
\date{}

\begin{document}
\maketitle

\begin{abstract}
In this paper we propose different notions of $\F_{\zeta}$-geometry,
for $\zeta$ a root of unity, generalizing notions of $\F_1$-geometry
(geometry over the ``field with one element") based on the behavior
of the counting functions of points over finite fields, the Grothendieck
class, and the notion of torification. We relate $\F_{\zeta}$-geometry
to formal roots of Tate motives, and to functions in the Habiro ring,
seen as counting functions of certain ind-varieties. We investigate
the existence of $\F_{\zeta}$-structures in examples arising from
general linear groups, matrix equations over finite fields, and some
quantum modular forms.
\end{abstract}

\tableofcontents

\section{Introduction}

This paper is a first step towards addressing a question asked by Don Zagier 
to the second author, about the possible existence of a notion of $\F_{\zeta}$-geometry,
with $\zeta$ a root of unity, generalizing the notion of $\F_1$-geometry (geometry over
the field with one element), and possibly related to his theory of quantum 
modular forms, \cite{Zagier}.

\subsection{The notion of $\F_1$-geometry}\label{F1Sec}

The idea of the ``field with one element" $\F_1$ arose from observations of
Jacques Tits on the behavior of the number of points of certain varieties defined over
finite fields $\F_q$, in the limit where $q\to 1$. For example, if $q=p^k$, with $p$ prime
and $\F_q$ is the finite field of characteristic $p$ with $q$ elements, then one has
$$ \# \P^{n-1}(\F_q) = \frac{\# (\A^n 
(\F_q) - \{ 0 \})}{\# \bG_m(\F_q)} =
 \frac{q^n-1}{q-1} = [n]_q $$
$$ \# {\rm Gr}(n,j)(\F_q) = \# \{ \P^j(\F_q) \subset \P^n(\F_q) \}  $$
$$ = \frac{[n]_q !}{[j]_q ! [n-j]_q !}
= \binom{n}{j}_q $$
where
$$ [n]_q ! = [n]_q [n-1]_q \cdots [1]_q, \ \ \  [0]_q ! =1  $$
Tits observed that in the limit where $q\to 1$ these expression become, respectively,
the cardinalities of the sets
$$ \P^{n-1}(\F_1) := \text{ finite set of cardinality } n $$
$$ {\rm Gr}(n,j) (\F_1) := \text{ set of subsets of cardinality } j $$

This observation suggested that, although it is impossible to define a field with
one element (by definition a field has at least two elements $0\neq 1$), finite
geometries often behave as if such a field $\F_1$ existed.
Further work of Kapranov-Smirnov, Soul\'e, Manin, and many others in more
recent years, developed different versions of geometry over $\F_1$.  For a 
comparative survey of different existing notions of $\F_1$-geometry, we refer
the reader to \cite{LL2}.

\subsection{The extensions $\F_{1^n}$}

The observation of Tits shows that the analog of an $n$-dimensional vector space $V$ 
in $\F_1$-geometry is just a finite set of cardinality $n$, and the analog of the 
linear transformations $\GL_n$ is the permutation group $S_n$. 

Kapranov and Smirnov further argued that one can make sense of ``extensions"
$\F_{1^n}$ of $\F_1$, analogous to the unramified extensions $F_{q^n}$ of the
finite fields $\F_q$. Their notion of a vector space over $\F_{1^n}$ is a pointed 
set $(V,v)$ with a free action of $\mu_n$  (the group of $n$-th roots of unity)
on $V- \{ v \}$, and linear maps become permutations compatible 
with the $\mu_n$-action, with the formal change of coefficients from 
$\F_1$ to $\Z$ given by
\[\mathbb{F}_{1^n} \otimes _{\mathbb{F}_1} \mathbb{Z}:=\mathbb{Z}[t,t^{-1}]/(t^n-1).\]

\subsection{Torifications and $\F_{1^n}$-geometry}

Among the various existing notions of $\F_1$ geometry, we focus here on
the one based on {\em torifications}, introduced by Lopez-Pe\~{n}a and Lorscheid
in \cite{LoPe}. The notion of torification can be regarded as a geometrization of
a decomposition of the Grothendieck class of the variety that reflects the existence
of $\F_1$ and $\F_{1^n}$ points, defined as limits of the counting of points over
$\F_{q^n}$.

\smallskip
\subsubsection{Counting points over $\F_{1^n}$}
Given a variety $X$ defined over $\Z$, one considers the
reductions $X_p$ modulo various primes $p$ and the corresponding counting
functions, for $q=p^r$, of the form $N_X(q)=\# X_p(\F_q)$. 

\smallskip

One obtains a good notion of $\F_{1^n}$-points
of $X$ when the variety satisfies the following conditions
(see Theorem 4.10 of \cite{CC} and Theorem 1 of \cite{Deit}).

\begin{cond}\label{FqF1cond}
Let $X$ be a variety over $\Z$. Then $X$ has $\F_{1^n}$-points for all $n\geq 0$ if
the following conditions hold:
\begin{enumerate}
\item {\em Polynomial countability:} The counting function $N_X(q)$ is a polynomial in $q$.
\item {\em Positivity:} The polynomial $N_X$ takes non-negative values $N_X(n+1)\geq 0$
for all $n\geq 0$.
\end{enumerate}
\end{cond}

The number of $\F_{1^n}$-points is given by 
\begin{equation}\label{F1npts}
\# X(\F_{1^n}) = N_X(n+1).
\end{equation}
The number of $\F_1$-points corresponds to the case $\# X(\F_1)=N_X(1)=\lim_{q\to 1} N_X(q)$.

\smallskip
\subsubsection{Grothendieck classes and torifications}
The polynomial countability condition is very closely related to
the nature of the motive of the variety $X$. In fact, one can lift
the two conditions above from the behavior of the counting of
points to the classes in the Grothendieck ring of varieties. 

\smallskip

The Grothendieck ring $K_0(\cV_\Z)$ of varieties over $\Z$ is
generated by isomorphism classes $[X]$ of varieties
with the relation
\begin{equation}\label{inclexcl}
[X] = [Y] + [X\smallsetminus Y],
\end{equation}
for any closed subvariety $Y\subset X$ and with the
product structure given by $$[X\times Y]=[X]\cdot [Y].$$
Classes in the Grothendieck ring are also called {\em virtual
motives} or {\em naive motives}.

\smallskip

It is customary to denote by $\bL=[\A^1]$ the Lefschetz motive,
the class of the affine line and by $\bT=\bL-1=[\bG_m]$ the
torus motive, the class of the multiplicative group. 

\smallskip

One can then strengthen the condition listed above for the
existence of $\F_{1^n}$-points in terms of the following
conditions on the Grothendieck class of a variety $X$ over $\Z$.

\begin{cond}\label{KXcond}
Let $X$ be a variety over $\Z$. The $X$ has a {\em motivic $\F_1$-structure} if
the following conditions are satisfied:
\begin{enumerate}
\item {\em Torification:} The class $[X]$ has a decomposition
\begin{equation}\label{torusKclass}
[X]=\sum_{k=0}^N a_k \, \bT^k.
\end{equation}
\item {\em Positivity:} All the coefficients $a_k$ of \eqref{torusKclass}
are non-negative integers, $a_k\geq 0$.
\end{enumerate}
\end{cond}

We refer to these conditions as a torification of the Grothendieck
class (see \cite{ManMar} for a discussion of the different levels of
torified structures). We refer to it as ``motivic $\F_1$-structure" as
it implies that the virtual motive $[X]$ of the variety behaves as
if the variety had an $\F_1$-structure in the stronger sense described
in Condition \ref{Torcond} below.

\smallskip

\begin{lem}\label{Cond21}
Condition \ref{KXcond} implies Condition \ref{FqF1cond}. 
\end{lem}

\proof As an additive
invariant, the counting of points factors through the Grothendieck
ring, hence we have
$$ N_X(q) = \sum_{k=0}^N a_k \, \# \bT^k(\F_q)=  \sum_{k=0}^N a_k \, (q-1)^k, $$
since $\# \bG_m(q)=q-1$. Moreover, the positivity $a_k\geq 0$ implies
$N_X(n+1)=\sum_k a_k \, n^k \geq 0$. 
\endproof

The number of $\F_1$-points
has a simple geometric interpretation in this setting as $X(\F_1)=a_0 = \chi(X)$,
the Euler characteristic of $X$. The relation between $\F_{1^n}$-points and roots
of unity is also evident, as the value $N_X(n+1)=\sum_k a_k \, n^k$  counts every
class of a zero-dimensional torus (there are $a_0$ of those, the $\F_1$-points) 
together with $n$ points in each $\bG_m$ (or $n^k$ points in each $\bG_m^k$)
which corresponds to counting the points of the group $\mu_n$ of roots of units 
of order $n$ embedded in the torus $\bG_m$.

\smallskip
\subsubsection{Geometric torifications}

A further strengthening of this condition is achieved through the notion of
a {\em geometric torification}. This is the notion of torifications considered in \cite{LoPe}.
One replaces the condition on the Grothendieck class with a geometric condition
on the variety, namely the requirement that the variety decomposes as a union of tori.

\begin{cond}\label{Torcond}
Let $X$ be a variety over $\Z$. The {\em geometric torification condition} is the
requirement that there is a morphism $e_X: T \to X$, for $T=\amalg_{j\in I} T_j$ 
with $T_j=\bG_m^{k_j}$, such that $e_X|_{T_j}$ is an immersion for each $j$ and 
$e_X$ induces bijections of sets of $k$-points, $T(k) \simeq X(k)$, over any field $k$.
\end{cond}

We also refer to this condition as a {\em geometric $\F_1$-structure}.

\smallskip

\begin{lem}\label{Cond32}
Condition \ref{Torcond} implies Condition \ref{KXcond}, hence also Condition \ref{FqF1cond}.
\end{lem}

\proof In fact, if $X$ has a geometric
torification in this sense, by the relation \eqref{inclexcl} in the Grothendieck ring, the
class decomposes as $[X]=\sum_{j\in I} \bT^{k_j}=\sum_k a_k \bT^k$, where the
coefficients satisfy $a_k\geq 0$, since $a_k =\#\{ T_j \subset T\,|\, T_j=\bG_m^k \}$. 
\endproof

\smallskip

One can impose further conditions on the torification, like the requirement that
it is {\em affine} or {\em regular}. We will not discuss those here, nor their specific
role in $\F_1$-geometry. We refer the reader to \cite{LoPe} and \cite{ManMar} for
further discussions of these properties and of the notion of morphisms of torified 
spaces. For our purposes here, we take a minimalistic approach to $\F_1$-geometry
and we only require the existence of a geometric torification as a condition for
the existence of an $\F_1$-structure on $X$.

\subsection{Towards a geometry over $\F_{\zeta}$}

Our purpose here is to use counting functions,
classes in the Grothendieck ring, and geometric torifications, to 
provide a possible setting for $\F_{\zeta}$-geometry, with $\zeta$ a root
of unity, compatible with the notions discussed above of $\F_1$ and 
$\F_{1^n}$-geometry.

\section{Two naive approaches to points over $\F_{\zeta}$}\label{NaiveSec}

We describe here two different approaches to a notion of points over $\F_{\zeta}$,
where $\zeta$ is a root of unity of order $n$.
The first is based on the idea of replacing the notion of $\F_1$-points given by
the limiting behavior when $q\to 1$ of the counting function $N_X(q)$ of $\F_q$-points
with the behavior of the same function at roots of unity. The second notion is instead
based on extending the polynomial interpolation of the counting function $N_X(q)$,
whose values at positive integers counts the $\F_{1^n}$-points to values at negative
integers as counting of $\F_{\zeta}$-points. 

\smallskip

In order to distinguish different notions of $\F_\zeta$-points, we will
be assigning different names to the various conditions we introduce.
So we will be talking of ``evaluation $\F_\zeta$-points" in the first case
and of ``interpolation $\F_\zeta$-points" in the latter, so that one can
keep in mind what procedure we are following.

\smallskip
\subsection{$\F_\zeta$-points and cell decompositions}

\begin{cond}\label{Fzetacond}
Let $\zeta$ be a root of unity, $\zeta =\exp(2\pi i k/n)$, 
and let $X$ be a variety over $\Z$, with
$N_X(q)$ the counting function of points over $\F_q$.
We say that $X$ has ``evaluation $\F_\zeta$-points" if the following conditions hold.
\begin{enumerate}
\item {\em Polynomial countability:} The counting function $N_X(q)$ is a polynomial 
$N_X(q)=\sum_{k=0}^N b_k\, q^k$.
\item {\em Positivity and divisibility:} the only coefficients $b_k\neq 0$ occur 
at degrees $k$ with $n|k$, with $b_k>0$.
\end{enumerate}
\end{cond}

Notice that the second condition implies a {\em positivity} condition for the
number of $\F_\zeta$-points, defined as the value $N_X(\zeta)$. In fact,
Condition \ref{Fzetacond} implies that the polynomial $N_X(q)$ 
has non-negative value $N_X(\zeta)\geq 0$.

\smallskip

The value $N_X(\zeta)$ will then count the number of $\F_\zeta$-points. 

As in the case of $\F_1$ and $\F_{1^n}$-points, we can make
this condition more geometric, by first rephrasing it in terms of
Grothendieck classes and then in terms of geometric decompositions
of the variety.

\begin{cond}\label{KXzetacond}
Let $\zeta$ be a root of unity, $\zeta =\exp(2\pi i k/n)$.
A variety $X$ over $\Z$ has an {\em evaluation $\F_\zeta$-structure} at the
{\em motivic level} if the following conditions hold:
\begin{enumerate}
\item {\em Affine decomposition:} The class $[X]$ has a decomposition
\begin{equation}\label{Affdec}
[X]=\sum_{k=0}^N b_k\, \bL^k
\end{equation}
\item {\em Positivity and divisibility:} the only coefficients $b_k\neq 0$ occur 
at degrees $k$ with $n|k$, with $b_k>0$.
\end{enumerate}
\end{cond}

At a more geometric level we can formulate a stronger geometric version of this
condition as a decomposition of the variety which implies the decomposition of
the Grothendieck class.

\begin{cond}\label{Affcond}
Let $\zeta$ be a root of unity, $\zeta =\exp(2\pi i k/n)$.
A variety $X$ over $\Z$ has an {\em evaluation $\F_\zeta$-structure} at
the {\em geometric level} if
there exists a cell decomposition $X =\amalg_{j\in I} \bA^{k_j}$
where the dimensions of the cells satisfy $n|k_j$ for all $j\in I$.
\end{cond}

Since we want to think of $\F_1$ as being in some sense an extension of $\F_\zeta$
we need to check that conditions we use to construct $\F_\zeta$-points and
$\F_\zeta$-structures imply the corresponding $\F_1$-conditions.

\begin{lem}\label{FzF1cond}
Condition \ref{Fzetacond} implies Condition \ref{FqF1cond}. Condition \ref{KXzetacond}
implies Condition \ref{KXcond} and Condition \ref{Fzetacond}. 
Condition \ref{Affcond} implies Condition \ref{KXzetacond} and Condition \ref{Torcond}.
\end{lem}

\proof The first condition is the same. The second condition
implies that $N_X(q)=\sum_{k=0}^N b_k\, q^k$ with non-negative
coefficients, hence $N_X(n+1)\geq 0$ for all $n\geq 0$. In the same
way, using $\bT=\bL-1$, one sees that the existence of a decomposition
\eqref{Affdec} into powers of the Lefschetz motive with non-negative
coefficients implies the existence of a decomposition \eqref{torusKclass} with
non-negative coefficients, since each $\bL^k = \sum_j \binom{k}{j} \bT^j$.
At the geometric level, the existence of a cell decomposition implies
the existence of a torification by decomposing each affine space
into a union of tori. The implications between the three $\F_\zeta$-conditions
follow as in the corresponding implications of the $\F_1$-conditions.
\endproof


\begin{rem}\label{remfinzeta} {\rm
With this notion, a variety $X$ can have ``evaluation $\F_\zeta$-points" for at most
finitely many $\zeta$'s, since the order of $\zeta$ has to be a number 
that divides all the exponents of the polynomial $N_X(q)$. }
\end{rem}

As we will see later, 
it is convenient to also consider a weaker form of Condition \ref{Fzetacond} and of
the corresponding geometrizations, Condition \ref{KXzetacond}
and Condition \ref{Affcond}, defined as follows.

\begin{cond}\label{Fzetaweakconds}
Let $\zeta$ be a root of unity, $\zeta =\exp(2\pi i k/n)$.
A variety $X$ over $\Z$ has a {\em partial evaluation $\F_\zeta$-structure} 
if the following holds:
\begin{itemize}
\item (Counting Function): $N_X(q)=\sum_{k=0}^N b_k\, q^{nk} + (q^n-1)\,P(q)$, with
$b_k\geq 0$ and a polynomial $P(q)$ satisfying $P(n+1)\geq 0$ for all $n\geq 0$.
\item (Grothendieck Class): $[X]=\sum_{k=0}^N b_k\, \bL^{nk} + (\bL^n-1) P(\bL)$,
with $b_k\geq 0$ and a polynomial $P(\bL)=\sum_j a_j \bT^j$ with $a_j\geq 0$.
\item (Geometric Decomposition): There is a subvariety $Y\subset X$ with
a cell decomposition $Y=\cup_{j\in J} \A^{nk_j}$, such that the
complement $X\smallsetminus Y$ admits a geometric torification,
with Grothendieck class $[X\smallsetminus Y]=(\bL^n-1) P(\bL)$.
\end{itemize}
\end{cond}

\begin{rem}\label{constrtor} {\rm The existence of a cell decomposition
into affine spaces implies the existence of a geometric torification.
Requiring that the complement of a torified submanifolds 
inside a torified manifold also has a torification is in general a very strong 
condition, see the discussion of ``complemented subspaces" and
constructible torifications in \cite{ManMar}. The condition on the existence
of a torification on $X\smallsetminus Y$ can be weakened by requiring
the existence of a {\em constructible torification} in the sense of \cite{ManMar},
with Grothendieck class $[X\smallsetminus Y]=(\bL^n-1) P(\bL)$. }
\end{rem}

\medskip
\subsection{$\F_{\zeta}$-points and polynomial interpolation}\label{FzetainterpSec}

One arrives at a different possible notion of $\F_{\zeta}$-points
if one extends the polynomial interpolation $N_X(n+1)=\# X(\F_{1^n})$
used to define points over the extensions $\F_{1^n}$ in such a way
that $\F_1$ would play the same role of an extension of $\F_{\zeta}$.
This leads to regarding the values of the counting function $N_X(q)$
at negative integers as representing points over $\F_{\zeta}$.

\begin{cond}\label{FzF1cond}
Let $X$ be a variety over $\Z$. Then $X$ has ``interpolation $\F_{\zeta}$-points" if
the following holds:
\begin{enumerate}
\item {\em Polynomial countability:} The counting function $N_X(q)$ is a polynomial in $q$.
\item {\em Positivity:} The polynomial $N_X$ takes non-negative values $N_X(m+1)\geq 0$
for all $m\geq 0$ and for $m=-n<0$.
\end{enumerate}
\end{cond}

\begin{rem}\label{infzetarem}
With this notion, a variety $X$ has can have ``interpolation $\F_{\zeta}$-points" for all $n$
if the polynomial $N_X$ has non-negative values at all (positive and negative)
integers.
\end{rem}

\smallskip

This notion of $\F_{\zeta}$-points can also be strengthened to
a condition on the Grothendieck class and on geometric decompositions
of the variety. We focus here on varieties that have an interpolation $\F_{\zeta}$-structure
simultaneously for all $n$. We formulate geometric conditions that extend the
notion of torification used in the $\F_1$-case. 

\begin{cond}\label{KXzetacond2}
A variety $X$ over $\Z$ has an {\em interpolation $\F_{\zeta}$-structure} at the
{\em motivic level}, for all $n\geq 0$, if
the following conditions hold:
\begin{enumerate}
\item {\em Torification:} The class $[X]$ has a decomposition
\begin{equation}\label{Tor1dec}
[X]=\sum_{k=0}^N a_k\, \bT^k
\end{equation}
\item {\em Dual torification:} The polynomial $P(\bL)=\sum_k a_k \, (-1)^k (\bL+1)^k$, obtained
by replacing $\bL \mapsto -\bL$ in \eqref{Tor1dec} can also be written in the form
\begin{equation}\label{dualTor1dec}
 P(\bL)=\sum_k c_k \, \bT^k, 
\end{equation} 
with coefficients $c_k\geq 0$.
\end{enumerate}
\end{cond}

\begin{ex}\label{exTdual}
Any polynomial $\sum_k b_k \bL^k$ that is even, with $b_k\geq 0$ and
$b_{2k+2}\geq b_{2k+1}$, satisfies \eqref{Tor1dec} and \eqref{dualTor1dec}.
\end{ex}

Condition \ref{KXzetacond2} implies the existence of interpolation $\F_\zeta$-points for
all $n$.

\begin{lem}\label{KX2FxF1}
Condition \ref{KXzetacond2} implies that Condition \ref{FzF1cond} is satisfied
for all $n\geq 0$.
\end{lem}

\proof Condition \eqref{Tor1dec} implies that $N_X(n+1)\geq 0$ for all $n\geq 0$.
By Condition \eqref{dualTor1dec}, 
the polynomial $\hat N_X(q):=N_X(-q)$ is then of the form $\hat N_X(q)=\sum_k c_k (q-1)^k$,
with all $c_k\geq 0$, hence positivity is satisfied at all negative integers as well, $\hat N_X(n)=
N_X(-n)\geq 0$, for all $n\geq 0$.
\endproof

This suggests a geometric condition that induces Condition \ref{KXzetacond2}
for the classes.

\begin{cond}\label{DualTorcond}
A variety $X$ over $\Z$ has an {\em interpolation $\F_{\zeta}$-structure} 
at the {\em geometric level}, for all $n\geq 0$, if
$X$ has a geometric torification $X=\cup_{j\in J} \bG_m^{k_j}$ and there
exists another variety $\hat X$ with a geometric torification 
$\hat X=\cup_{\ell \in I} \bG_m^{d_\ell}$, such that the Grothendieck classes
are related by
\begin{equation}\label{dualKclass}
 [\hat X]= \sum_k (-1)^k b_k \, \bL^k, \ \ \ \text{ where } \ \ \  [X]=\sum_k b_k \bL^k. 
\end{equation} 
\end{cond}

\begin{ex}\label{P2Ndual}
For $X=\P^{2N}$ a variety satisfying \eqref{dualKclass}
is $\hat X= \A^0 \cup \bigcup_{k=1}^N (\A^{2k}\smallsetminus \A^{2k-1})=
\A^0 \cup \bG_m \times (\bigcup_{k=1}^N \A^{2k-1}) $. 
\end{ex}

The different approaches described in this section, naively
based on the ideas of evaluation at $q=\zeta$ and of interpolation
over negative integers of the positivity condition $N_X(n+1)\geq 0$,
seem to lead to very different notions of $\F_\zeta$-geometry.
In the following sections we describe a motivic setting which
will allow us to reconcile and unify these views.

\section{The Habiro ring and the Grothendieck ring}\label{HabSec}

A more sophisticated way to keep track of conditions,
at the level of classes in the Grothendieck ring, that
would correspond to ``$n$-th roots of $\F_1$-structures"
can be achieved by considering cyclotomic completions
of the Grothendieck ring of varieties, modeled on the
well known Habiro completion of polynomial rings, \cite{Habiro},
\cite{Habiro2}, \cite{Manin}. This type of completion allows
for a consistent system of evaluations and expansions
at all roots of unity. We show that, in our motivic setting,
the target of the evaluation maps are Grothendieck
rings of a family of orbit categories of pure motives,
with respect to the tensor action of the Tate motives $\Q(n)$.

\subsection{The Habiro ring as analytic functions on roots of unity}

The Habiro ring was introduced in \cite{Habiro} as a {\em cyclotomic completion}
of the polynomial ring $\Z[q]$, or more generally of $R[q]$ with $R$ a commutative
ring with unit. Namely, one defines the Habiro ring $\widehat{\Z[q]}$ 
(also sometimes denoted $\Z[q]^\N$) as the projective limit 
\begin{equation}\label{HabRing}
\widehat{\Z[q]}= \varprojlim_n \Z[q] / ((q)_n)
\end{equation}
of the quotients $\Z[q] / ((q)_n)$ by the principal
ideals $\cI_n=((q)_n)$ generated by the elements
$(q)_n:=((1-q)(1-q^2)\cdots (1-q^n))$. The ideals $\cI_n$ are ordered by 
divisibility, with $(q)_k|(q)_n$ for $k\leq n$, so that $((q)_n)\subset ((q)_k)$,
with projections $\Z[q]/((q)_n)\twoheadrightarrow \Z[q]/((q)_k)$.

\smallskip
\subsubsection{Formal series}

The elements of $\widehat{\Z[q]}$ are therefore formal series $\sum_n P_n(q)$
with $P_n(q) \in \cI_n$, or equivalently formal series
\begin{equation}\label{eltsHab}
a_0(q) + \sum_{n\geq 1} a_n(q) \, (1-q)(1-q^2)\cdots (1-q^n),
\end{equation}
with the $a_n(q) \in \Z[q]$ for $n\geq 0$.

\smallskip

As observed in \cite{Manin}, one can interpret the Habiro ring as a deformation
of the ring $\hat\Z$, given by the projective limit (see \cite{Dantz}, \cite{Len})
\begin{equation}\label{hatZunivnum}
\hat\Z = \varprojlim_n \Z/(n!),
\end{equation}
with the product $(q^n-1)\cdots (q-1)$ replacing $n!$.

\smallskip
\subsubsection{Evaluation at roots of unity}

Given any roots of unity $\zeta$, there is an evaluation map, which
is a surjective ring homomorphism
\begin{equation}\label{evzeta}
e_\zeta: \widehat{\Z[q]} \to \Z[\zeta].
\end{equation}
Moreover, as shown in Theorem 6.2 of \cite{Habiro}, the
map obtained by considering all evaluations, for all roots of unity,
\begin{equation}\label{evmap}
e: \widehat{\Z[q]} \to \prod_\zeta \Z[\zeta], \ \ \ \ e(f) = (e_\zeta(f))_\zeta=(f(\zeta))_\zeta
\end{equation}
is {\em injective}. In fact, a more refined result in \cite{Habiro} shows that
evaluations on certain infinite sets of roots of unity suffices to determine
the function in the Habiro ring. 

\smallskip
\subsubsection{Taylor expansions at roots of unity}

Functions in the Habiro ring have
well defined Taylor expansions at all roots of unity, given by ring homomorphisms
\begin{equation}\label{Taylor}
\sigma_\zeta: \widehat{\Z[q]} \to \Z[\zeta] [[q-\zeta]],
\end{equation}
and the function in the Habiro ring is completely determined 
by its Taylor expansion at any root of unity, since by Theorem
5.2 of \cite{Habiro} these homomorphisms are injective.
These properties can be used to interpret the Habiro ring
as a ring of {\em analytic functions at roots of unity}.
More recently, Manin interpreted the Habiro ring as a possible
approach to analytic geometry over $\F_1$, \cite{Manin}. In a
similar perspective, the Habiro ring was related to endomotives
and the Bost--Connes system in \cite{Mar}.

\subsection{Projective systems of Grothendieck rings}

For our purposes, we will consider here a construction similar
to the Habiro ring, but based on the Grothendieck ring of varieties.

\begin{defn}\label{HabGrdef}
Let $(\bL)_n:=(\bL-1)(\bL^2-1)\cdots (\bL^n-1)$, namely the class
$$ (\bL)_n =[ (\A^n\smallsetminus \{ {\bf 0} \})\times \cdots \times (\A^2\smallsetminus \{ {\bf 0} \})\times
(\A^1\smallsetminus \{ 0 \}) ], $$
where ${\bf 0}=(0,\ldots,0) \in \A^k$.
Letting $\cI_n=((\bL)_n)$ be the ideal in $K_0(\cV_\Z)$ generated
by this element, we obtain a quotient
$$ K_0(\cV_\Z)/\cI_n = K_0(\cV_\Z)/((\bL-1)(\bL^2-1)\cdots (\bL^n-1)). $$
Since $((\bL)_n) \subset ((\bL)_k)$, for $k\leq n$, the ideals $\cI_n$
are ordered by inclusion, with projections $K_0(\cV_\Z)/\cI_n\twoheadrightarrow
K_0(\cV_\Z)/\cI_k$ and we can consider the projective limit
\begin{equation}\label{projlimK0}
\widehat{K_0(\cV_\Z)} := \varprojlim_n K_0(\cV_\Z)/((\bL-1)(\bL^2-1)\cdots (\bL^n-1)). 
\end{equation}
We refer to this ring as the {\em Habiro--Grothendieck ring}. 
\end{defn}

Elements in this ring
can be written as formal series
\begin{equation}\label{eltsHabGr}
\alpha_0 + \sum_{n\geq 1} \alpha_n \, (\bL-1)(\bL^2-1)\cdots (\bL^n-1), \ \ \  \alpha_k \in K_0(\cV_\Z),
\, k\geq 0.
\end{equation}

\smallskip
\subsubsection{Habiro--Grothendieck ring and the Tate motive}

In the theory of motives it is customary to introduce a formal inverse of the Lefschetz
motive, the Tate motive $\bL^{-1}=:\Q(1)$. Correspondingly, one considers the Grothendieck
ring $K_0(\cV_\Z)[\bL^{-1}]$. In our setting, we can then consider the Habiro--Grothendieck
ring given by
\begin{equation}\label{2projlimK0}
\widehat{K_0(\cV_\Z)[\bL^{-1}]} = \varprojlim_n K_0(\cV_\Z)/((\bL-1)(\bL^2-1)\cdots (\bL^n-1)). 
\end{equation}

\begin{lem}\label{HabK0mot}
Let $[\GL_n] \in K_0(\cV_\Z)$ be the Grothendieck class of the general linear 
group $\GL_n$. The ring $\widehat{K_0(\cV_\Z)[\bL^{-1}]}$ of \eqref{2projlimK0} 
is equivalently obtained in the forms
\begin{equation}\label{HabmotGLn}
\widehat{K_0(\cV_\Z)[\bL^{-1}]}=\varprojlim_n K_0(\cV_\Z)[\bL^{-1}]/([\GL_n]),
\end{equation}
\end{lem}

\proof We have (see for instance Lemma 2.6 of \cite{Bridge})
\begin{equation}\label{GLnclass}
 [\GL_n] = \bL^{n(n-1)/2} (\bL-1)(\bL^2-1)\cdots (\bL^n-1). 
\end{equation} 
\endproof

\smallskip
\subsubsection{The Tate motive as a Habiro--Grothendieck class}
The versions \eqref{projlimK0} and \eqref{2projlimK0}, \eqref{HabmotGLn}
of the Habiro--Grothendieck ring are in fact equivalent, since the Lefschetz
motive $\bL$ already has an inverse $\bL^{-1}$ in \eqref{projlimK0}.

\begin{lem}\label{TateHabGr}
In the Habiro--Grothendieck ring $\widehat{K_0(\cV)}$ the Lefschetz
motive $\bL$ is invertible, with the Tate motive $\bL^{-1}$ given by
$$ \bL^{-1} = \sum_{n\geq 0} \bL^n (\bL^n-1)\cdots (\bL^2-1)(\bL-1). $$
\end{lem}

\proof The argument is as in Proposition 7.1 of \cite{Habiro}.
\endproof

\smallskip

\begin{rem}\label{GLquotinv} {\rm 
If instead of taking quotients with respect to the
ideals generated by the classes $[\GL_n]$, one inverts these
elements, one obtains $K_0(\cV_\Z)[ [\GL_n]^{-1}, n\geq 1]$,
which is isomorphic to the Grothendieck ring of algebraic stacks, \cite{Bridge}.}
\end{rem}

\smallskip
\subsubsection{The Habiro ring of Tate motives}

The subring of virtual Tate motives in $K_0(\cV_\Z)$ (respectively, in
$K_0(\cV_\Z)[\bL^{-1}]$) is the polynomial ring $\Z[\bL]$ (respectively,
the ring of Laurent polynomials $\Z[\bL,\bL^{-1}]$). The corresponding
Tate Habiro--Grothendieck ring is the classical Habiro ring 
\begin{equation}\label{TateHabGr}
\widehat{\Z[\bL]}=\varprojlim_n \Z[\bL]/((\bL-1)(\bL^2-1)\cdots (\bL^n-1)),
\end{equation}
whose elements are formal series
\begin{equation}\label{eltsHabGr}
f_0(\bL) + \sum_{n\geq 1} f_n(\bL) \, (\bL-1)(\bL^2-1)\cdots (\bL^n-1), \ \ \  
f_k \in \Z[\bL], \, k\geq 0.
\end{equation}

\smallskip
\subsubsection{Evaluation maps}

For $n\geq 1$, let $(\bL^n-1)$ be the ideal in $K_0(\cV_\Z)$
generated by the element $\bL^n-1=[\A^n \smallsetminus \{ 0 \}]$.
This determines a quotient $K_0(\cV_\Z)/(\bL^n-1)$.  The inclusion
$((\bL)_n)\subset (\bL^n-1)$ implies the existence of surjective
evaluation maps as follows.

\begin{lem}\label{evK0}
There are evaluation maps, in the form of surjective ring
homomorphisms
\begin{equation}\label{evK0maps}
ev_n: \widehat{K_0(\cV_\Z)} \twoheadrightarrow K_0(\cV_\Z)/(\bL^n-1).
\end{equation}
\end{lem}

We give a better motivic interpretation of the target $K_0(\cV_\Z)/(\bL^n-1)$
of the evaluation map in terms of orbit categories of motives.

\subsection{Orbit categories of motives}

Given an additive category $\cC$ and an automorphism $F$, the {\em orbit category}
has ${\rm Obj}(\cC/F)={\rm Obj}(\cC)$ and morphisms given by
\begin{equation}\label{morOrbCat}
{\rm Hom}_{\cC/F}(X,Y):= \oplus_{k\in \Z} {\rm Hom}_{\cC}(X, F^k(Y)).
\end{equation}
In particular, if $\cC$ is also a symmetric monoidal category and $F=-\otimes \cO$
with $\cO$ a $\otimes$-invertible object in $\cC$, the orbit category $\cC/_{-\otimes\cO}$
also has a symmetric monoidal structure compatible with the projection functor, see
\S 7 of \cite{Tab}.

\smallskip
\subsubsection{Tate motives and orbit categories}

In particular, we consider here the case where $\cM={\rm Chow}_\Q$, the
category of Chow motives with rational coefficients and $\cO=\Q(n)$,
where $\Q(1)=\bL^{-1}$ is the Tate motive and $\Q(n)=\Q(1)^{\otimes n}$.
The case with $n=1$ gives rise to the orbit category $\cM_1={\rm Chow}_\Q/_{-\otimes \Q(1)}$
considered in \cite{Tab}, which embeds in the category of noncommutative motives.
Here we also consider the orbit categories $\cM_n={\rm Chow}_\Q/_{-\otimes \Q(n)}$ with
$n>1$. 

\smallskip
\subsubsection{Grothendieck groups of orbit categories}

The Grothendieck group $K_0(\cC)$ of a pseudo-abelian category $\cC$ is the quotient of
the free abelian group on the isomorphism classes $[M]$ of objects in $\cC$
by the subgroup generated by elements $[M]-[M']-[M'']$ for $M\simeq M'\oplus M''$ in $\cC$.
In the symmetric monoidal case, $K_0(\cC)$ is also a ring. 

\begin{prop}\label{K0orbnMot}
Let $K_0(\cM)$ be the Grothendieck ring of the category $\cM={\rm Chow}_\Q$ of Chow
motives. Then the Grothendieck ring of the orbit category $\cM_n={\rm Chow}_\Q/_{-\otimes \Q(n)}$
can be identified with the quotient 
\begin{equation}\label{K0Mn}
K_0(\cM_n)=K_0(\cM)/((\bL^n-1)).
\end{equation}
\end{prop}

\proof Isomorphic objects $M\simeq N$ in ${\rm Chow}_\Q/_{-\otimes \Q(n)}$
are objects $M,N$ in ${\rm Chow}_\Q$ related by an isomorphism
in ${\rm Hom}_{{\rm Chow}_\Q/_{-\otimes \Q(n)}}(M,N)$, that is, such that 
$M \simeq N \otimes \Q(nk)$ in ${\rm Chow}_\Q$, for some $k\in \Z$.
This introduces new relations of the form $[M]=[M] \bL^{nk}$ in the Grothendieck
ring: these are elements of the ideal $(\bL^n -1)$. All relations 
$[M]-[M']-[M'']$ in $K_0({\rm Chow}_\Q)$ are still relations
in $K_0({\rm Chow}_\Q/_{-\otimes \Q(n)})$. Thus, we can identify
$K_0({\rm Chow}_\Q/_{-\otimes \Q(n)})$ with the quotient 
$K_0({\rm Chow}_\Q)/(\bL^n-1)$.
\endproof

\smallskip
\subsubsection{Euler characteristic and evaluation}
There are motivic Euler characteristics associated to the orbit categories 
${\rm Chow}_\Q/_{-\otimes \Q(n)}$.

\begin{cor}\label{motEulch}
The motivic Euler characteristic $\chi_{mot}$ of \cite{GiSou} induces
motivic Euler characteristics $\chi_{mot}^{(n)}(X)$ of varieties $X$,
which determine ring homomorphisms $\chi_{mot}^{(n)}:
K_0(\cV)/((\bL^n-1)) \to K_0(\cM)/((\bL^n-1))$.
\end{cor}

\proof
Gillet and Soul\'e constructed a motivic Euler characteristic $\chi_{mot}$, which
is a ring homomorphism $\chi_{mot}: K_0(\cV) \to K_0(\cM)$, by associating
to a variety $X$ a complex $W(X)$ in ${\rm Chow}_\Q$, whose class $[W(X)]\in 
K_0(\cM)$ defines $\chi_{mot}(X)$. Further mapping $W(X)$ to its class
$[W(X)]\in K_0({\rm Chow}_\Q/_{-\otimes \Q(n)})$ defines a motivic Euler
characteristic $\chi_{mot}^{(n)}: K_0(\cV) \to K_0(\cM)/((\bL^n-1))$. 
The homomorphism $\chi_{mot}$ maps $\bL \in K_0(\cV)$
to $\bL\in K_0(\cM)$, and the induced morphism descends to a ring 
homomorphism $\chi_{mot}^{(n)}: K_0(\cV)/((\bL^n-1)) 
\to K_0(\cM)/((\bL^n-1))$.
\endproof

This, the evaluation maps \eqref{evK0maps} can be seen as
maps to the Grothendieck rings of all the orbit categories 
$\cM_n ={\rm Chow}_\Q/_{-\otimes \Q(n)}$. We will see in \S \ref{TaterootSec}
below that some natural operations on the Habiro ring described in
\cite{Mar}, modelled on the Bost--Connes quantum statistical mechanical
system, have a motivic interpretation in terms of ``roots of Tate motives". 

\medskip

\subsection{Counting functions and the Habiro ring}

Over a finite field $\F_q$, the counting function $X \mapsto N_q(X)=\# X(\F_q)$
factors through the Grothendieck ring $N_q: K_0(\cV_{\F_q})\to \Z$, since it satisfies
the inclusion-exclusion and product relations. 
The counting function induces an identification 
$N: \widehat{\Z[\bL]} \to \widehat{\Z[q]}$ of the Tate
part $\widehat{\Z[\bL]}\subset \widehat{K_0(\cV_\Z)}$ of 
the Habiro--Grothendieck ring with the Habiro ring on
polynomials $\Z[q]$, and we can interpret functions $f(q)$
in $\widehat{\Z[q]}$ as counting functions of Tate classes
in $\widehat{K_0(\cV_\Z)}$.

\smallskip

Under the evaluation maps \eqref{evK0maps}, the counting functions 
are correspondingly mapped to the a counting function in $\Z[q]/(q^n-1)$
as the image of the counting function $\widehat{\Z[q]}$, consistently with
the observation in \cite{Tab}, that counting functions on the orbit category
of motives ${\rm Chow}_\Q/_{-\otimes \Q(1)}$ are only defined up to $q-1$.

\section{Roots of Tate motives}\label{TaterootSec}

The idea of introducing Tate motives $\Q(r)$ for $r$ not an integer
was suggested by Manin in \cite{Manin2}. In particular, in the
case of finite fields, a nice geometric interpretation of $\Q(1/2)$, 
a square root of the Tate motive $\Q(1)$, in terms of a 
supersingular elliptic curves was suggested by Manin
and fully developed by Ramachandran in \cite{Rama}. Other
considerations on the ``exotic" Tate motives $\Q(r)$ with $r\notin \Z$,
can be found in \S 7 of \cite{Den}. Here we will take a purely formal
approach to the construction of these objects, motivated by 
a semigroup of endomorphisms of the Habiro ring studied in \cite{Mar}
from a quantum statistical mechanical point of view.

\subsection{Roots of Tate motives and the Habiro--Grothendieck ring}

We focus here on the Tate part $\Z[\bL]$ of the Grothendieck ring
of varieties $K_0(\cV)$ and the Habiro ring $\widehat{\Z[\bL]}$.

As in \cite{Mar}, we consider endomorphisms 
$\sigma_n: \widehat{\Z[\bL]} \to \widehat{\Z[\bL]}$
of the Habiro ring induced by the morphisms $\sigma_n: \Z[\bL]\to \Z[\bL]$
given by $\bL \mapsto \bL^n$. These extend to the Habiro
completion since $(\bL)_m | \sigma_n(\bL)_m$, see Proposition 2.1
of \cite{Mar}. 

\smallskip

As in \cite{Mar}, we denote by $\widehat{\Z[\bL]}_\infty$ the direct
limit 
\begin{equation}\label{dirlimsigman}
\widehat{\Z[\bL]}_\infty = \varinjlim ( \sigma_n: \widehat{\Z[\bL]} \to \widehat{\Z[\bL]} )
\end{equation}
of the system of maps $\sigma_n: \widehat{\Z[\bL]} \to \widehat{\Z[\bL]}$.
As shown in Lemma 2.3 and Proposition 2.2 of \cite{Mar}, the ring $\widehat{\Z[\bL]}_\infty$
has an equivalent description as follows.

\smallskip

Let $\cP_\Z=\Z[\bL^r; r\in \Q^*_+]$ the ring of polynomials in the rational 
powers of the variable $\bL$, and let $\widehat{\cP_\Z}$ be the completion
\begin{equation}\label{hatPZ}
\widehat{\cP_\Z} = \varprojlim_N \cP_\Z /\cJ_N,
\end{equation}
with respect to the ideals $\cJ_N$ generated by the $(\bL^r)_N=(\bL^{rN}-1)\cdots (\bL^r-1)$
with $r\in \Q^*_+$. Then the direct limit above can be identified with
\begin{equation}\label{hatZinftyPZ}
\widehat{\Z[\bL]}_\infty \simeq \widehat{\cP_\Z}.
\end{equation}
The morphisms $\sigma_n: \widehat{\Z[\bL]} \to \widehat{\Z[\bL]}$ induce automorphisms
of $\widehat{\Z[\bL]}_\infty$, which induce an action of the group $\Q^*_+$ by
\begin{equation}\label{Qact}
f(\bL) \mapsto  f(\bL^r), \ \ \  r\in \Q^*_+.
\end{equation}
This action can be encoded in the crossed product ring $\widehat{\Z[\bL]}_\infty\rtimes \Q^*_+$
considered in \cite{Mar}.

\smallskip

In the motivic setting we are considering, this means that taking the
direct limit of the Tate part $\widehat{\Z[\bL]}$ of the Habiro--Grothendieck
ring with respect to the endomorphisms $\sigma_n: \widehat{\Z[\bL]} \to \widehat{\Z[\bL]}$
has the effect of introducing formal roots $\bL^r$ of the Lefschetz motive, with $r\in \Q^*_+$.
We can view this at the level of categories of motives, by introducing formal
roots of Tate motives $\Q(r)$, with $r\in \Q^*_+$ using the Tannakian formalism.
We will discuss the relevance of roots of Tate motives to $\F_\zeta$-geometry
in the following section.

\subsection{Formal roots of Tate motives}

Let $\cT = {\rm Num}_\Q^\dagger$ denote the Tannakian category of pure motives,
with the numerical equivalence relation on cycles, with $\dagger$ denoting the
twist in the commutativity morphisms of the tensor structure by the Koszul sign 
rule, which is needed to make the category Tannakian. The category $\cT$ is
equivalent to the category of finite dimensional linear representation of an
affine group scheme, the motivic Galois group $G={\rm Gal}(\cT)$. The inclusion
of the Tate motives in ${\rm Num}_\Q^\dagger$ determines a group homomorphism
$t: G \to \bG_m$, such that $t\circ w=2$, where $w: \bG_m \to G$ is the weight
homomorphism (see \S 5 of \cite{DelMil}).

\smallskip

We proceed as in the construction of the square root of $\Q(1)$ given in \S 3.4
of \cite{KoSo}. This describes the category obtained by formally adding to $\cT$
a new object $\Q(1/p)$ with the property that $\Q(1/p)^{\otimes p}=\Q(1)$, namely
a formal tensor $p$-th root of the Tate motive.

\begin{defn}\label{Qproot}
For $\cT = {\rm Num}_\Q^\dagger$ and $G={\rm Gal}(\cT)$,
consider the homomorphisms $\sigma_p: \bG_m \to \bG_m$
given by $\sigma_p: \lambda \mapsto \lambda^p$, and let 
$G^{(p)}$ denote the fibered product 
\begin{equation}\label{Gp}
G^{(p)}=\{ (g,\lambda)\in G\times \bG_m\,:\, t(g)=\sigma_p(\lambda) \}.
\end{equation} 
The resulting Tannakian category ${\rm Rep}_{G^{(p)}}$ is denoted by
$\cT(\Q(\frac{1}{p}))$.
\end{defn}

We can iterate this construction over the set of generators of the
multiplicative semigroup $\N$ given by the set of primes.

\begin{lem}\label{Gpq}
Let $p\neq q$ be prime numbers and consider the homomorphisms 
$\sigma_p$ and $\sigma_q$ of $\bG_m$ as above. Let $G^{(p,q)}$ and
$G^{(q,p)}$ denote, respectively, the fiber products 
$$ G^{(p,q)}=\{ (g,\lambda) \in G^{(q)}\times \bG_m\,:\, t(\pi_q(g))=\sigma_p(\lambda) \}, $$ 
$$ G^{(q,p)}=\{ (g,\lambda) \in G^{(p)}\times \bG_m\,:\, t(\pi_p(g))=\sigma_q(\lambda) \}, $$
with $\pi_p: G^{(p)}\to G$ and $\pi_q:G^{(q)}\to G$ the projections of the fibered product.
Then $G^{(p,q)}\simeq G^{(q,p)}$, with ${\rm Rep}(G^{(p,q)})$ the Tannakian category
obtained by adjoining formal roots $\Q(1/p)$ and $\Q(1/q)$ with $\Q(1/p)^{\otimes p}= \Q(1)
=\Q(1/q)^{\otimes q}$, and $\Q(1/(pq))$ with $\Q(1/(pq))^{\otimes p}=\Q(1/q)$
and $\Q(1/(pq))^{\otimes q}=\Q(1/p)$.
\end{lem}

\proof We can identify 
$$ G^{(p,q)}\simeq G^{(q,p)}\simeq \{ (g,\lambda,\tilde\lambda)\in G\times \bG_m\times \bG_m \,|\, 
t(g) =\sigma_p(\lambda)=\sigma_q(\tilde\lambda) \}, $$
hence the resulting pullback only depends on the unordered set $\{ p, q\}$. 
By construction, the pullback introduces additional generators in the 
category ${\rm Rep}(G^{(p,q)})$, of the form $\Q(1/p)$, $\Q(1/q)$,  $\Q(1/(pq))$,
with the listed properties under tensor powers.
\endproof

One can iterate this process, by considering the set of all primes in $\Z$,
the generators of the multiplicative semigroup $\N$ is positive integers,
and, associate to the first $N$ elements $\{ p_1,\ldots, p_N \}$ the
corresponding group $G^{(p_1,\ldots, p_N)}$ obtained as the iterated
pullback along the map $\sigma_{p_j}: \bG_m\to \bG_m$ of
$G^{(p_1,\ldots, p_{j-1})}$, for $j=1,\ldots,N$. The corresponding
Tannakian category $\cT^{(p_1,\ldots, p_{j-1})}$ is obtained from
$\cT={\rm Num}_\Q^\dagger$ by introducing additional generators
$\Q(1/n)$ for all $n\in \N$ with primary decomposition $n=p_1^{a_1}\cdots p_N^{a_N}$
and with the obvious relations under tensor products,
$$ \Q(\frac{1}{p_1^{a_1}\cdots p_N^{a_N}})^{\otimes p_j}=
\Q(\frac{1}{p_1^{a_1}\cdots p_j^{a_j-1}\cdots p_N^{a_N}}). $$
The following result is then a direct consequence of the construction.

\begin{cor}\label{rootTatecat}
One can then consider the limit $G^\N =\varprojlim G^{(p_1,\ldots, p_{j-1})}$.
The corresponding category of representations $\cT^\N$ contains roots of Tate
motives $\Q(1/n)$ of all order $n\in \N$. 
The subcategory of $\cT^\N$ generated
by all the roots $\Q(1/n)$ of the Tate motive contains the usual Tate motives
and has an action of the group $\Q^*_+$ by automorphisms. 
\end{cor}

\medskip
\subsection{Roots of Tate motives and orbit categories}

There are surjective evaluation maps from the direct limit $\widehat{\Z[\bL]}_\infty$ of
\eqref{dirlimsigman} and the quotient rings $\Z[\bL^r]/(\bL^{Nr}-1)$, for any given
$r\in \Q^*_+$ and $N\in \N$. In particular, it suffices to look at the cases where
$r=1/p$ and $(N,p)=1$. The quotient ring $\Z[\bL^{1/p}]/(\bL^{N/p}-1)$ can be seen as
the $K_0$ of the Tate part of the orbit category $\cT(\Q(1/p))/_{-\otimes \Q(N/p)}$,
by the same argument as in Proposition \ref{K0orbnMot}.
Namely, we have the following, which follows directly from the previous discussion.

\begin{prop}\label{K0orbTQp}
The category $\cT(\Q(1/p))$ has
$$ K_0(\cT(\Q(1/p)))=K_0(\cT)[t]/(t^q-\bL), $$
and the subcategory of Tate motives, generated by $\bL^{1/p}$, has $K_0$ given by
$\Z[\bL^{1/p}]$. The orbit category $\cT(\Q(1/p))/_{-\otimes \Q(N/p)}$ has
$$ K_0(\cT(\Q(1/p))/_{-\otimes \Q(N/p)})=K_0(\cT(\Q(1/p)))/(\bL^{N/p}-1) $$
and the subcategory of Tate motives has $K_0$ given by
$\Z[\bL^{1/p}]/(\bL^{N/p}-1)$. There are surjective homomorphisms
$$ ev_{p,N} : \widehat{\Z[\bL]}_\infty \twoheadrightarrow K_0(\cT(\Q(1/p))/_{-\otimes \Q(N/p)}). $$
\end{prop}

\section{Roots of Tate motives and $\F_\zeta$-geometry}\label{RootsFzetaSec}

We can now introduce another possible way of thinking of $\F_\zeta$ structures,
which will be related to the notion of ``evaluation $\F_\zeta$-structure" discussed
in \S \ref{NaiveSec}.

\begin{cond}\label{FzetaQroots}
Let $\zeta$ be a root of unity of order $n$. 
Let $X$ be a variety over $\Z$ that has an $\F_1$-structure, with Grothendieck class
$[X]=\sum_k b_k \,\bL^k$ in $\Z[\bL]$. Then a {\em Tate root $F_\zeta$-structure} 
for $X$ is a choice of an object $M$ in the category $\cT^\N$ of motives containing all
roots of Tate motives, with the property that the class $[M]\in K_0(\cT^\N)$ is of the form
$$ [M] =\sum_k b_k \, \bL^{k/n}. $$
\end{cond}

We can view the class $[M]$ as being in $K_0(\cT(\Q(1/n))$ rather than in
the larger category $\cT^\N$. 

\smallskip

\begin{rem}\label{evTateFzeta}{\rm 
The condition on the existence of an ``evaluation $\F_\zeta$-structure" discussed in
\S \ref{NaiveSec} above becomes the condition that there exists a 
``Tate root $F_\zeta$-structure" $M$, such that $M$ is in fact in $\cT$ itself,
inside the larger $\cT(\Q(1/n))$. }
\end{rem}

\smallskip

\begin{ex}\label{affcellTateFzeta}{\rm If $X$ has an affine cell decomposition
$X=\amalg_{j\in J} \A^{k_j}$, then $[X]=\sum_j \bL^{k_j}$ and an example of a 
choice of a Tate root $F_\zeta$-structure is given by $M=\oplus_{j\in J} \bL^{k_j/n}$.}
\end{ex}

\section{$\F_\zeta$-geometry in the Habiro--Grothendieck ring}\label{FzetaHabSec}

We now upgrade our previous notions of $\F_\zeta$ points to
versions that can be formulated for functions in the Habiro ring,
viewed as counting functions for classes in the Tate part of
the Grothendieck--Habiro ring, and then further upgrade the
resulting conditions to geometric conditions based on
torification, for certain classes of ind-varieties. 

\smallskip
\subsection{Ind-varieties and Habiro functions}

Recall that an ind-variety over $\Z$ is a direct limit $X=\varinjlim_\alpha X_\alpha$ of a
direct system of varieties $X_\alpha$ over $\Z$ with morphisms 
$\phi_{\alpha,\beta}: X_\alpha \to X_\beta$ for all $\alpha<\beta$ in the directed set of
indices $\alpha,\beta\in \cI$.

\smallskip

Suppose given a formal series 
$f(q)=\sum_{m=0}^\infty a_m(q) (q^m-1)\cdots (q^2-1) (q-1)$
in the Habiro ring $\widehat{\Z[q]}$, with the property that 
the polynomials $a_m(q)$ are counting functions
\begin{equation}\label{amqcount}
a_m(q)= N_{X_m}(q),
\end{equation}
for some varieties $X_m$ over $\Z$. 
Consider then, for each $N\geq 0$ the varieties $\cX_N$ given by
\begin{equation}\label{cXNind}
\cX_N = \bigcup_{m=0}^{N-1} X_m \times (\A^m\smallsetminus \{ {\bf 0} \})\times
\cdots (\A^2 \smallsetminus \{ {\bf 0} \})\times (\A^1\smallsetminus \{ 0 \}).
\end{equation}
This gives formally an interpretation (not necessarily unique) of the Habiro function
as a counting function of an ind-variety $f(q)=N_{\cX}(q)$. We can use this interpretation
to formulate conditions on the existence of $\F_\zeta$ structures based on properties of
functions in the Habiro ring.

\smallskip
\subsection{$\F_1$ structures on ind-varieties}

We first recast the notion of $\F_1$-structure based on torification
and the behavior of the counting function, for ind-varieties and 
classes in the Habiro--Grothendieck ring. We formulate the condition 
from the strongest geometric level to the weaker level of the counting 
function in the Habiro ring.

\begin{cond}\label{F1Habindvar}
An ind-variety $\cX=\varprojlim_N \cX_N$ over $\Z$ has an 
$\F_1$-structure if the following holds:
\begin{itemize}
\item (Geometric Decompositions): The varieties $\cX_N$ have
a decomposition of the form \eqref{cXNind}, where the varieties $X_m$ 
admit a geometric torification, $X_m=\cup_{j\in J} \bT^{k_j}$.
\item  (Grothendieck Class):  The class $[\cX_N]$ defines an element
in the Tate part of the Habiro--Grothendieck ring $\widehat{\Z[\bL]}$,
of the form $\sum_{m=0}^\infty \alpha_m(\bL) (\bL^m-1)\cdots (\bL-1)$,
where $[X_m]=\alpha_m(\bL)$ is a polynomial in $\Z[\bL]$ of the form
$\alpha_m(\bL)=\sum_{j\in J} \bT^{k_j}$.
\item (Counting Function): The counting functions $N_{\cX_N}(q)$
determine a function $N_{\cX}(q)$ in the Habiro ring, with
$N_{\cX}(q)=\sum_{m=0}^\infty \alpha_m(q) (q^m-1)\cdots (q-1)$,
where $\alpha_m(q)=N_{X_m}(q)=\sum_{j\in J} (q-1)^{k_j}$.
\end{itemize}
\end{cond}

In particular, the counting functions $N_{X_m}(q)$ and $N_{\cX_N}(q)$
have the property that the values at all the non-negative integers $n\geq 0$
are non-negative, and provide a counting of points over all the extensions
$\F_{1^n}$.

\smallskip
\subsection{$\F_\zeta$ structures on ind-varieties}

We now refine Condition \ref{F1Habindvar} for the existence of an $\F_1$-structure
to a notion of $\F_\zeta$ structure based on our previous notion of ``interpolation"
$\F_\zeta$-points, described in \S \ref{FzetainterpSec} above.

\begin{cond}\label{FzetaHabindvar}
An ind-variety $\cX=\varprojlim_N \cX_N$ over $\Z$ has an 
$\F_\zeta$-structure, for $\zeta^n=1$, if it has an $\F_1$-structure
in the sense of Condition \ref{F1Habindvar} and the polynomial
\begin{equation}\label{evnpoly}
N_{\cX_n}(q)=\sum_{m=0}^{n-1} \alpha_m(q) (q^m-1)\cdots (q-1)
\end{equation}
satisfies $N_{\cX_n}(-n) \geq 0$.
In particular, an ind-variety $\cX$ has an $\F_\zeta$-structure for
roots of unity $\zeta$ of all orders if $N_{\cX_n}(-n) \geq 0$ for all $n\geq 0$.
\end{cond}

The counting function $N_{\cX_n}(q)$ modulo $q^n-1$ 
gives the counting function associated to the class
$ev_n([\cX]) \in K_0(\cM)/(\bL^n-1)$, under the evaluation
maps \eqref{evK0maps}.

\smallskip
\subsection{Constructible torifications and constructible $\F_\zeta$-structures}\label{ConstrSec}

In \cite{ManMar} a weaker form of the torification condition of \cite{LoPe} was
introduced: the notion of a {\em constructible torification}, which defines a structure
of $\F_1$-constructible set, instead of the stronger notion of $\F_1$-variety.
As in ordinary algebraic geometry, where constructible sets are complements
of algebraic varieties inside other varieties (which are not always varieties themselves),
so the notion of a ``constructible" torification encodes the operation of taking 
complements of geometric torifications inside other torifications, where the
complement is not necessarily itself torified. We recall here the main definitions.

Recall first that, as defined in \cite{ManMar}, one considers the set $\cC_{\F_1}$
of all constructible sets over $\Z$ that can be obtained, starting with the
multiplicative group $\bG_m$ by repeated operations of products, disjoint
unions, and taking complements. 

\begin{defn}\label{constor}
A variety (or more generally a 
constructible set ) $X$ defined over $\Z$ has a {\em constructible torification} if
the following conditions hold:
\begin{enumerate}
\item There is a morphism of constructible sets 
$e: C \to X$ from a $C\in \cC_{\F_1}$ such that the restriction of $e$ to each
component of $C$ is an immersion and the map $e(k): C(k)\to X(k)$ 
is a bijection of the sets of $k$-points for any field $k$. 
\item The class $[X]\in K_0(\cV_\Z)$ has a decomposition $[X]=\sum_k a_k \bT^k$
with $a_k\geq 0$.
\end{enumerate}
\end{defn}

We reformulate this condition in the context of ind-varieties 
and Habiro functions, and we obtain the following.

\begin{cond}\label{constorind}
An ind-variety $\cX=\varprojlim_N \cX_N$ over $\Z$ has a
{\em constructible} $\F_1$-structure if the following holds:
\begin{itemize}
\item (Geometric Decompositions): The varieties $\cX_N$ have
a decomposition of the form \eqref{cXNind}, where the varieties $X_m$ 
admit a {\em constructible torification} as in Definition \ref{constor}.
\item  (Grothendieck Class):  The class $[\cX_N]$ defines an element
in the Tate part of the Habiro--Grothendieck ring $\widehat{\Z[\bL]}$,
of the form $\sum_{m=0}^\infty \alpha_m(\bL) (\bL^m-1)\cdots (\bL-1)$,
where $[X_m]=\alpha_m(\bL)$ is a polynomial in $\Z[\bL]$ of the form
$\alpha_m(\bL)=\sum_k a_{m,k} \bT^k$
with $a_{m,k}\geq 0$.
\item (Counting Function): The counting functions $N_{\cX_N}(q)$
determine a function $N_{\cX}(q)$ in the Habiro ring, with
$N_{\cX}(q)=\sum_{m=0}^\infty \alpha_m(q) (q^m-1)\cdots (q-1)$,
where $\alpha_m(q)=N_{X_m}(q)=\sum_k  a_{m,k} (q-1)^k$,
with $a_{m,k}\geq 0$.
\end{itemize}
\end{cond}

\medskip
\subsection{Roots of Tate motives and $\F_\zeta$ structures in the
Habiro--Grothendieck ring}

As in \S \ref{RootsFzetaSec}, we can also consider the notion
of ``Tate root $\F_\zeta$-structure" in the setting of the Habiro--Grothendieck ring.
For a class 
$$ [\cX] = \sum_{m=1}^\infty \alpha_m(\bL) (\bL^m-1)\cdots (\bL-1) $$
a Tate root $\F_\zeta$-structure, for $\zeta$ a root of unity of order $n$, is an element 
$$ f(\bL^{1/n})=\sum_{m=1}^\infty \alpha_m(\bL^{1/n}) (\bL^{m/n}-1)\cdots (\bL^{1/n}-1) $$
in the ring $\widehat{\Z[\bL]}_\infty=\hat\cP_\Z$ of \eqref{dirlimsigman}.

\section{Examples}

In this section we discuss several explicit examples of functions defining
$\F_\zeta$-structures on associated ind-varieties in the sense described above.

Some of the examples considered in this section are quantum modular forms
in the sense of Zagier \cite{Zagier}, but we do not know at this stage whether
there is an interpretation of the quantum modularity condition in terms of some
of the notions of $\F_\zeta$-geometry that we have been considering in this
paper. We hope to return to this question elsewhere.

\subsection{General linear groups}

It is well known that, as in \eqref{GLnclass},
\begin{equation}\label{GLnpts}
\# \GL_m(\F_q) = q^{m(m-1)/2}(q-1)(q^2-1)\dots (q^m-1).
\end{equation}
We consider a function in the Habiro ring, associated to the
general linear groups, of the form
\begin{equation}\label{GLfq}
 f_{GL}(q):= 1+ \sum_{m=1}^{\infty}q^{m(m-1)/2}(q-1)(q^2-1)\dots (q^m-1) .
\end{equation}
This is the counting function for the class in the Grothendieck--Habiro ring
\begin{equation}\label{GLGrHab}
[ \cup_m \GL_m ] =\sum_{m=0}^\infty [\GL_m],
\end{equation}
where we set $[\GL_0]=1$.

\smallskip

\begin{prop}\label{GLFzeta}
Let $\cX_n =\amalg_{m=0}^n \GL_m$. The counting function $N_{\cX_n}(1-n)$
is non-negative for odd $n$, hence $\cX=\amalg_{m=0}^\infty \GL_m$ has
an $\F_\zeta$ structure for all roots of unity of odd order.
\end{prop}

\proof We consider the sum as starting at $m=1$.
Summing the $m$th and $(m+1)$th terms, we have 
$q^{m(m-1)/2}(q-1)(q^2-1)\dots (q^m-1)[q^m(q^{m+1}-1)+1]$.
When $m=4k$, we have $(q-1)(q^2-1)\dots (q^m-1) |_{q=1-n}
\ge(-2)^{2k}\times1^{2k}\ge 2$. We also have 
$q^{m(m-1)/2}|_{q=1-n}=q^{2k(4k-1)}|_{q=1-n}\ge(-2)^{2k(4k-1)}\ge 2$,
while $(1-n)^{m+1}-1\le -3$ and $(1-n)^m\ge 2$, so $(q^m(q^{m+1}-1)+1)|_{q=1-n}\le-5$,
hence the sum of these terms is negative. When $m=4k+1$, 
$((q-1)(q^2-1)\dots (q^m-1))|_{q=1-n}$ is negative, since there is an odd number
of negative terms, while $(1-n)^{m(m-1)/2}=(1-n)^{2k(4k+1)}$ is positive, and
$(1-n)^{m+1}-1\ge1$ and $(1-n)^m<- 1$, so $(1-n)^m((1-n)^{m+1}-1)+1<0$. Thus,
the sum of the terms is positive in this case. When $m=4k+3$,
$((q-1)(q^2-1)\dots (q^m-1))|_{q=1-n}$ is positive, since there is an 
even number of negative terms, while $(1-n)^{m(m-1)/2}=(1-n)^{(2k+1)(4k+3)}$
is negative and $(1-n)^{m+1}-1\ge1$ and $(1-n)^m<- 1$, so $(1-n)^m((1-n)^{m+1}-1)+1<0$.
Thus, the sum of the terms if positive in this case. We then look at the cases for $n$.
For $n=1$ the sum starting from $m=1$ is zero, and for $n=2$ it is negative.
For $n$ even, the sum of terms larger than $m=0$  begins with an odd term 
and ends with an odd term.  Since the sum of the $m$th and $(m+1)$th term 
is negative for $m$ even, the sum is bounded by $\sum_{m=2}^{n-1}N_{\cX_n}(1-n)
\le -20\times\frac{n-2}{2}$. Adding the terms $m=0$ and $m=1$ does not change
sign, hence the sum is negative. When $n$ is odd, the sum of terms larger than $m=0$ 
begins with an odd term and ends with an even term.  Since the sum of the $m$th 
and $(m+1)$th term is positive for $m$ odd, the sum is positive for odd $n$.
\endproof

\subsection{Matrix equations over $\F_q$}

As another example, we consider a special case of the matrix equations over
finite fields analyzed in \cite{Carl}. Let $A$ be a non-singular matrix of order $2m$.  
Then, as in \cite{Carl}, one knows 
that the number of solutions $Z(A,A)$ for $X'AX=A$ is given by
\begin{equation}\label{matrixsol}
E_{2m}(q)=Z(A,A)=q^{m^2}\prod_{i=1}^m(q^{2i}-1).
\end{equation}
The term with $m=0$ is set equal to $1$.
One sees immediately, that one can also read the expression \eqref{matrixsol}
as counting points of a variety of the following form.

\begin{lem}\label{E2mcount}
The number of solutions \eqref{matrixsol} is the counting function $N_{X_m}(q)$ of
the variety
$X_m=\mathbb{A}^{m^2}\times
(\mathbb{A}^2-\{ {\bf 0 }\})\times(\mathbb{A}^4-\{ {\bf 0 }\})
\times\dots\times(\mathbb{A}^{2m}-\{ {\bf 0 }\})$.
\end{lem}

We now consider the positivity condition $N_{\cX_n}(1-n)\geq 0$ for the existence
of an $\F_\zeta$-structure as in Condition \ref{FzetaHabindvar}. 

\begin{prop}\label{E2mposn}
The counting function $N_{\cX_n}(q)$ satisfies $N_{\cX_n}(1-n)\geq 0$ for odd $n$ and
for $n$ of the form $4k+2$. Hence the ind-variety $\cX$, in this case, has an $\F_\zeta$
structure for all roots of unity of odd order and of order $4k+2$.
\end{prop}

\proof Summing the $(m-1)$th term and the $m$th term for $m\ge2$, 
we have $(q^2-1)\dots(q^{2k-2}-1)(q^{m^2}(q^{2m}-1)-q^{(m-1)^2})$.
We need to check whether $(1-n)^{m^2}((1-n)^{2m}-1)+(1-n)^{(m-1)^2}$ 
is positive. When $m$ is odd, $m\geq 3$, we have $(1-n)^{(m-1)^2}>1$,
and the positivity of the above term is the same as the positivity of 
$[(1-n)^\text{odd}((1-n)^\text{even}-1)+1]\leq -2)\times(a\in\mathbb{Z}^+)+1\le-1$,
hence the sum of the two terms is negative. When $m$ is even, the sign of the
same term is given by the sign of $(1-n)^\text{odd}[(1-n)^\text{odd}((1-n)^\text{even}-1)+1]
\ge(-2)(-1)=2$, hence the sum if positive. Thus, when $n$ is odd, the sum
$1+\sum_{m=1}^{n-1} N_{X_m}(1-n)\ge1+(2)\times\frac{n-1}{2}\ge n>0$. When
$n=4k$, for $k\ge1$, one has $\sum_{m=2}^{n/2-1}N_{X_m}(1-n)
\le(-1)\times\frac{n/2-2}{2}=1-k$, hence the sum $\sum_{m=0}^{n/2-1}N_{X_m}(1-n)\leq
2-k+(1-n)(n^2-2n)\le -k-4<0$. When $n=4k+2$, the final term $\frac{n}{2}-1$ is even,
hence the sum behaves as in the odd case, ending with an even term, hence the
sum is again positive.
\endproof

\subsection{Ramanujan's $q$-hypergeometric function} 

We consider next the example of the two $q$-hypergeometric functions
\[ \sigma (q)=\sum_{n=0}^{\infty} \frac{q^{n(n+1)/2}}{(1+q)(1+q^2)\dots (1+q^n)}=1+q-q^2+2q^3-2q^4+q^5+q^7-2q^8+ \dots\]
\[ \sigma^* (q)=2\sum_{n=1}^{\infty} \frac{(-1)^nq^{n^2}}{(1-q)(1-q^3)\dots (1-q^{2n-1})} = -2q-2q^2-2q^3+2q^7+2q^8+2q^{10}+ \dots\]
The function $\sigma(q)$ can be interpreted as he counting function for a class
in the Habiro--Grothendieck ring in the following way.

\begin{lem}\label{qhyplem}
The $q$-hypergeometric function
\begin{equation}\label{qhyp}
 \sigma(q) = \sum_{n=0}^\infty \frac{q^{n(n+1)/2}}{(1+q)(1+q^2)\cdots (1+q^n)} 
\end{equation}
is the counting function of the  Habiro--Grothendieck class
\begin{equation}\label{sigmaHabGr}
\sigma(\bL) = \bL^0 + \sum_{n=0}^\infty \bL^{n+1} (\bL-1)(\bL^2-1)\cdots (\bL^n-1).
\end{equation}
There is an ind-variety 
\begin{equation}\label{Xvarsigma}
\cX_\sigma=\bigcup_{n=0}^\infty \, \A^{n+1}\times (\A^1\smallsetminus \{ 0 \})\times (\A^2\smallsetminus \{ 0 \})
\times \cdots \times (\A^n\smallsetminus \{ 0 \})
\end{equation}
such that $[\cX_\sigma]=\sigma(\bL)$. All the varieties 
$$ \cX_{\sigma,N} =  \bigcup_{n=0}^N \,
\A^{n+1}\times (\A^1\smallsetminus \{ 0 \})\times (\A^2\smallsetminus \{ 0 \})
\times \cdots \times (\A^n\smallsetminus \{ 0 \}) $$
admit a geometric torification.
\end{lem}

\proof The $q$-hypergeometric functions \eqref{qhyp}
satisfies the identity  (\cite{ADH}, see also \cite{Zagier})
\begin{equation}\label{sigmaq}
\sigma(q)=1 + \sum_{n=0}^\infty q^{n+1} (q-1)(q^2-1) \cdots (q^n-1).
\end{equation}
In the form \eqref{sigmaq} this function can be thought of as the
counting function associated to the class \eqref{sigmaHabGr}.
It is clear by construction that $[\cX_\sigma]=\sigma(\bL)$. To show
that the varieties $\cX_{\sigma,N}$ admit a geometric torification
it suffices to show that a complement $\A^k\smallsetminus \{ 0 \}$
admits a torification, as one can then construct torifications on 
all the products 
$\A^{n+1}\times (\A^1\smallsetminus \{ 0 \})\times (\A^2\smallsetminus \{ 0 \})
\times \cdots \times (\A^n\smallsetminus \{ 0 \})$. We proceed inductively.
For $k=1$ we have $\A^1\smallsetminus \{ 0 \}=\bG_m$ which is already a torus.
Let us assume that we have constructed a torification for the complement
$\A^k\smallsetminus \{ 0 \}$. Then 
we decompose $\A^{k+1}\smallsetminus \{ 0 \}$ as a
disjoint union of $(\A^{k+1} \smallsetminus \A^k)\cup (\A^k\smallsetminus \{ 0 \})$.
The first term can be identified with $\A^k \times \bG_m$ where both factors
admit torifications, while the second term has a torification by inductive hypothesis.
\endproof

\smallskip

\begin{prop}\label{XsigmaFzeta}
The counting functions $N_{\cX_{\sigma,n}}(q)$ satisfy $N_{\cX_{\sigma,n}}(1-n)\geq 0$
for $n=1$, $n =4k+3$, and $n=4k$. Thus $\cX_\sigma$ has an $\F_\zeta$ structure at 
nontrivial roots of unity $\zeta$ of order $4k+3$ and $4k$.
\end{prop}

\proof Summing the $m$th, $(m+1)$th, $(m+2)$th and $(m+3)$th terms
in the series
\[ \sigma (q)-1=\sum_{m=0}^{\infty}q^{m+1}(q-1)(q^2-1)\dots (q^m-1)\]
gives
$$(q-1)(q^2-1)\dots (q^m-1)(q^{m+4}(q^{m+3}-1)(q^{m+2}-1)(q^{m+1}-1) $$
$$ +q^{m+3}(q^{m+2}-1)(q^{m+1}-1)+q^{m+2}(q^{m+1}-1)+q^{m+1})= $$
$$ q^{m+1}(q-1)(q^2-1)\dots (q^m-1)(q^3(q^{m+3}-1)(q^{m+2}-1)(q^{m+1}-1) $$
$$ +q^2(q^{m+2}-1)(q^{m+1}-1)+q(q^{m+1}-1)+1).$$ 

When $m=4k$, we have
$$q^2(q^{m+2}-1)(q^{m+1}-1)+q(q^{m+1}-1)|_{q=1-n}= $$ 
$$ q(q^{m+1}-1)(q(q^{m+2}-1)+1)|_{q=1-n}\le
-q(q^{m+1}-1)|_{q=1-n}\le-6,$$ 
$$(q^3(q^{m+3}-1)(q^{m+2}-1)(q^{m+1}-1)+1)|_{q=1-n}
\le-18+1=-17. $$ 
So we have 
$$
(q^3(q^{m+3}-1)(q^{m+2}-1)(q^{m+1}-1)+q^2(q^{m+2}-1)(q^{m+1}-1)+q(q^{m+1}-1)+1)|_{q=1-n}\le-23. $$ 
The term $(q-1)(q^2-1)\dots (q^{m}-1)|_{q=1-n}\ge2$, since 
there is an even number of negative terms, and $(1-n)^{m+1} \le-2$. Thus, the sum of the four
terms is bounded below by $92$, hence positive. 

\smallskip

In the case $m=4k+1$, we have
$$ q^2(q^{m+3}-1)(q^{m+2}-1)+q(q^{m+2}-1)=q(q^{m+2}-1)(q(q^{m+3}-1)+1)$$ 
with 
$(q(q^{m+3}-1)+1)|_{q=1-n}\le-1$ and $$q^2(q^{m+3}-1)(q^{m+2}-1)+q(q^{m+2}-1)
 =(q(q^{m+2}-1)(q(q^{m+3}-1)+1))|_{q=1-n}\le-6,$$ $$(q^2(q^{m+3}-1)(q^{m+2}-1)
+q(q^{m+2}-1)+1)|_{q=1-n}\le-5,$$ 
$$(q^3(q^{m+3}-1)(q^{m+2}-1)(q^{m+1}-1)+q^2(q^{m+2}-1)(q^{m+1}-1)+q(q^{m+1}-1)=$$
$$ q(q^{m+1}-1)(q^2(q^{m+3}-1)(q^{m+2}-1)+q(q^{m+2}-1)+1))|_{q=1-n}\ge10, $$
$$
(q^3(q^{m+3}-1)(q^{m+2}-1)(q^{m+1}-1)+q^2(q^{m+2}-1)(q^{m+1}-1)+q(q^{m+1}-1)+1)|_{q=1-n}\ge11, $$ 
while $(q-1)(q^2-1)\dots (q^{m}-1)|_{q=1-n}\le-2$, since there is an odd number of negative terms, and $(1-n)^{m+1} \ge2$. So the sum of the four terms is bounded above by $-44$, hence  negative. 

\smallskip

For $m=4k+2$, we have $$(q^2(q^{m+2}-1)(q^{m+1}-1)+q(q^{m+1}-1)=q(q^{m+1}-1)
(q(q^{m+2}-1)+1))|_{q=1-n}\le-q(q^{m+1}-1)\le-6,$$  $$(q^3(q^{m+3}-1)(q^{m+2}-1)
(q^{m+1}-1)+1)|_{q=1-n}\le -17,$$ so that $$(q^3(q^{m+3}-1)(q^{m+2}-1)(q^{m+1}-1)+q^2(q^{m+2}-1)(q^{m+1}-1)+q(q^{m+1}-1)+1)|_{q=1-n} \le-23.$$ The product $(q-1)(q^2-1)\dots (q^{m}-1)|_{q=1-n}\le-2$ is positive since there is an even number of negative terms and $(1-n)^{m+1} \le-2$. Thus, the sum
of the four terms is bounded above by $-92$, hence negative.

\smallskip

For $m=4k+3$, in $q^2(q^{m+3}-1)(q^{m+2}-1)+q(q^{m+2}-1)=q(q^{m+2}-1)(q(q^{m+3}-1)+1)$
we have $(q(q^{m+3}-1)+1)|_{q=1-n} \le-1$ and $$(q^2(q^{m+3}-1)(q^{m+2}-1)+q(q^{m+2}-1))|_{q=1-n}=(q(q^{m+2}-1)(q(q^{m+3}-1)+1))|_{q=1-n}\le-6,$$ 
$$(q^2(q^{m+3}-1)(q^{m+2}-1)+q(q^{m+2}-1)+1)|_{q=1-n}\le-5,$$ 
$$(q^3(q^{m+3}-1)(q^{m+2}-1)(q^{m+1}-1)+q^2(q^{m+2}-1)(q^{m+1}-1)+q(q^{m+1}-1))|_{q=1-n}= $$
$$ (q(q^{m+1}-1)(q^2(q^{m+3}-1)(q^{m+2}-1)+q(q^{m+2}-1)+1))|_{q=1-n}\ge10,$$ 
$$(q^3(q^{m+3}-1)(q^{m+2}-1)(q^{m+1}-1)+q^2(q^{m+2}-1)(q^{m+1}-1)+q(q^{m+1}-1)+1)|_{q=1-n}\ge11.$$
The product $(q-1)(q^2-1)\dots (q^{m}-1)|_{q=1-n} \ge2$ since there is an even number of 
negative terms, and $(1-n)^{m+1} \ge2$. Thus, the sum of the four terms is bounded below by
$44$, hence positive.

\smallskip

When $n=1$, $\sigma(1-n)=2>0$. When $n=2$, $\sigma(1-n)=1+\sum_{m=0}^{1}q^{m+1}(q-1)(q^2-1)\dots (q^m-1)|_{q=1-n}=1+q+q^2(q-1)|_{q=1-n}=-2<0$.

\smallskip

When $n\ge3$ and of the form  $n=4k+1$, the sum from $m=1$ to $n-1=4k$ can be 
split into groups of four terms each, with each group beginning with $m=4k+1$, hence we obtain  
$\sigma(1-n)\le-44k+1<0$.

\smallskip

For $n=4k+2$, the sum from $m=2$ to $n-1=4k+1$ can be split into groups of four terms,
each starting with $m=4k+2$, hence $\sigma(1-n)\le-92k+1+(q+q^2(q-1))|_{q=1-n}
\le-92k-n-6<0$.

\smallskip

When $n=4k+3$, the sum from $m=3$ to $n-1=4k+2$ can be split into groups of 
four terms starting with $m=4k+3$, and one has $\sigma(1-n)\ge 44k+1+
(q+q^2(q-1)+q^3(q-1)(q^2-1))|_{q=1-n}$, where 
$(q^2(q-1)+q^3(q-1)(q^2-1))|_{q=1-n}=q^2(q-1)(q(q^2-1)+1)|_{q=1-n}\ge6$, so that
$\sigma(1-n)\ge 44k+1+(q+q^2(q-1)+q^3(q-1)(q^2-1))|_{q=1-n}\ge 44k+7+1-4k-3=40k+5>0$.

\smallskip

When $n=4k$, the sum from $m=4$ to $n-1=4k-1$ can be separated into 
groups of four terms each starting with $m=4k$, so that one obtains  
$\sigma(1-n)\ge 92k+1+(q+q^2(q-1)+q^3(q-1)(q^2-1)+q^4(q-1)(q^2-1)(q^3-1))|_{q=1-n}
\ge92k+25+1-4k=88k+26>0$.
\endproof

\smallskip

We also consider the case of the function $\sigma^*(q)$, which we write in the
equivalent form (see \cite{Zagier})
\[\sigma^*(q)=2\sum_{k=0}^{\infty}(-1)^{k+1} q^{k+1}(q^2-1)(q^4-1) \dots (q^{2k}-1). \]
We can view this as the counting function of some ind-variety with an $\F_1$-structure 
in the following way.

\begin{lem}\label{F1sigmastar}
Let $\cX_{\sigma^*,\ell}$ be a union of tori with Grothendieck class
$$ [\cX_{\sigma^*,\ell}] = 4\ell\, \bT + \sum_{k=2}^{4\ell+1} \binom{4\ell+1}{k} \bT^k, $$
with $\bT=\bL-1$. Then setting 
$$ \cY_\ell=  \A^{2\ell}\times (\A^2\smallsetminus \{ {\bf 0} \}) \times \cdots
\times ( \A^{4\ell-2} \smallsetminus \{ {\bf 0} \}) \times \cX_{\sigma^*,\ell} $$
and $\cY=\cup_{\ell=1}^\infty \cY_\ell$, gives an ind-variety formally satisfying 
$$ 2 [\cX]= 2\sum_{\ell=1}^\infty \bL^{2\ell+1} (\bL^2-1) \cdots (\bL^{4\ell}-1) \, [\cX_{\sigma^*,\ell}]
=\sigma^*(\bL). $$
\end{lem}

\proof The difference of two consecutive terms in $\sigma^*(q)$ is
$$  q^{2\ell+1}(q^2-1)(q^4-1) \cdots (q^{4\ell}-1) - 
  q^{2\ell}(q^2-1)(q^4-1) \cdots (q^{4\ell-2}-1) = $$ $$ q^{2\ell} (q^2-1)\cdots (q^{4\ell-2}-1)  
  (q^{4\ell+1}-q-1) . $$
This can be seen as the counting function $N_{\cY_\ell}(q)$, since 
$\bL^{4\ell+1}-\bL-1 = 4\ell\, \bT + \sum_{k=2}^{4\ell+1} \binom{4\ell+1}{k} \bT^k$.
The varieties
$\A^{2\ell}$ and $(\A^k\smallsetminus \{ {\bf 0} \})$ have a geometric torification,
and so does  $\cX_{\sigma^*,\ell}$, as a union of tori. Thus, the varieties $\cY_\ell$
can be torified and $\cY$ has an $\F_1$-structure according to our definition.
\endproof

\begin{prop}\label{sigmastarFzeta}
The counting function $N_{\cY_\ell}(q)$ satisfied $N_{\cY_\ell}(1-n)\geq 0$ for all
$n=4k$, hence $\cY$ has an $\F_\zeta$ structure at all roots of unity $\zeta$
of order $4k$.
\end{prop}

\proof As above, summing the $(k-1)$th term and the $k$th term for $k\ge2$, 
we have $(q^2-1)\dots(q^{2k-2}-1)(q^{k+1}(q^{2k}-1)+q^{k})$. Thus, we 
need to check when $((1-n)^{k+1}((1-n)^{2k}-1)+(1-n)^{k})=(1-n)^k((1-n)
((1-n)^{2k}-1)+1)$ is positive. For all $k$, $(1-n)((1-n)^{2k}-1)+1)\le(-2)(1+1)\le-2$.
Moreover, for $k$ odd, $(1-n)^k\le-2$, hence $(1-n)^k[(1-n)((1-n)^{2k}-1)+1]\ge4$,
while when $k$ is even, $(1-n)^k\ge2$ and $(1-n)^k((1-n)((1-n)^{2k}-1)+1)\le-4$.

When $n=1$, $\sigma^*(q)|_{q=1-n}=0$, and when 
$n=2$, we have $\sigma^*(q)|_{q=1-n}=-2<0$.

For $n>2$ odd, the sum for $k$ from $1$ to $n-1$ begins with an odd term and ends with 
an even term, hence the sum is bounded above by $-4\times\frac{n-1}{2}$, and
$$ 2 \sum_{k=0}^{n-1} (-1)^{k+1} q^{k+1}(q^2-1)(q^4-1) \dots (q^{2k}-1) |_{q=1-n} 
\le-4(n-1)+2(1-n)=-6(n-1)<0. $$

For $n>2$ even, the sum for $k$ from $1$ to $n/2-1$ begins with an 
odd term and ends with $\frac{n}{2}-1$. If $n$ is of the form $4k+2$, then the
sum ends with an even term and the resulting summation 
$$ 2 \sum_{k=0}^{n/2-1} (-1)^{k+1} q^{k+1}(q^2-1)(q^4-1) \dots (q^{2k}-1) |_{q=1-n} $$
behaves as in the case $n$ odd discussed above. If $n$ is of the form $4k$ for $k\ge1$,
on the other hand, the sum ends in an odd term, and
$$ 2 \sum_{k=0}^{n/2-1} (-1)^{k+1} q^{k+1}(q^2-1)(q^4-1) \dots (q^{2k}-1) |_{q=1-n} \geq $$
$$ 8(k-1)+(2q+2q^2(q^2-1))|_{q=1-n} \geq 8k-4>0. $$
\endproof

\subsection{Kontsevich's function}

Finally, we show an example where the weaker notion of {\em constructible torification} 
occurs, as in \cite{ManMar}. The formal series
\begin{equation}\label{Kontser}
K(q)= \sum_{k=0}^\infty (1-q)(1-q^2)\cdots (1-q^k)
\end{equation}
was introduced by Kontsevich in relation to Feynman
integrals, and further studied by Zagier in \cite{Zagier2}.
It does not converge on any open set in $\bC$, but it has a well
defined value at all roots of unity. 

\smallskip

In our setting, \eqref{Kontser} can be interpreted as a counting
function, for an ind-variety that has a {\em constructible} torification.

\begin{prop}\label{Kontcount}
The formal series \eqref{Kontser} is the counting function for
the class in the Habiro--Grothendieck ring
\begin{equation}\label{KontClass}
K(\bL)= \sum_{k=0}^\infty (\bL^{2k}-2) (\bL^{2k-1}-1)\cdots (\bL^2-1)(\bL-1).
\end{equation}
We use the notation ${\bf 0}=(0,0,\ldots, 0)\in \A^k$ and
${\bf 1}=(1,1,\ldots, 1)\in \A^k$. The ind-variety 
\begin{equation}\label{XKont}
X_K = \cup_{k=0}^\infty (\A^{2k}\smallsetminus \{ {\bf 0}, {\bf 1} \})\times
(\A^{2k-1} \smallsetminus \{ {\bf 0} \}) \times \cdots 
(\A^2 \smallsetminus \{ {\bf 0} \}) \times (\A^1\smallsetminus \{ 0 \})
\end{equation}
has $[X_K]=K(\bL)$ and all the varieties $X_{K,N}$ given by
$$ X_{K,N} = \cup_{k=0}^N (\A^{2k}\smallsetminus \{ {\bf 0}, {\bf 1} \})\times
(\A^{2k-1} \smallsetminus \{ {\bf 0} \}) \times \cdots 
(\A^2 \smallsetminus \{ {\bf 0} \}) \times (\A^1\smallsetminus \{ 0 \}) $$
admit {\em constructible} torifications, in the sense of Section \ref{ConstrSec} above,
when $N=4k$.
\end{prop}

\proof We interpret the sum of two consecutive terms in the
Kontsevich series as 
$$ (1-q)(1-q^2)\cdots (1-q^{2k-1}) + (1-q)(1-q^2)\cdots (1-q^{2k})= $$
$$ (q^{2k}-2)(q^{2k-1}-1)\cdots (q-1), $$
which is the counting function for the product
$(\bL^{2k} -2)\times (\bL^{2k-1}-1) \times \cdots \times (\bL-1)$.
Clearly, by construction $K(\bL)$ is the class of $X_K$, and 
the existence of a constructible torification on $X_{K,N}$ can be shown
by showing that the complement $\A^{2k}\smallsetminus \{ {\bf 0}, {\bf 1} \}$
admits a constructible torification, as then a constructible torification of $X_{K,N}$
can be obtained using that of $\A^{2k}\smallsetminus \{ {\bf 0}, {\bf 1} \}$, together
with geometric torifications on the $\A^k\smallsetminus \{ {\bf 0} \}$, constructed as in
Lemma \ref{qhyplem}. The variety $\A^{2k}\smallsetminus \{ {\bf 0}, {\bf 1} \}$ clearly
admits a constructible torification, by taking the complement of a point $\{ {\bf 1} \}$ 
(with its trivial torification) in a torification $\A^{2k}\smallsetminus \{ {\bf 0} \}$. However,
the positivity of the Grothendieck class, which is necessary for the notion of
constructible torification is {\em not} satisfied by the individual term 
$\A^{2k}\smallsetminus \{ {\bf 0}, {\bf 1} \}$, since its Grothendieck class in the
variable $\bT$ is given by $-1 +\sum_{j=1}^{2k} \binom{2k}{j} \bT^j$. However,
the varieties $X_{K,N=4k}$ satisfy the positivity of the Grothendieck class, since
the negative terms cancel in pairs, hence the $X_{K,4k}$ have a 
constructible torification in the full sense of Definition \ref{constor}.
\endproof

\bigskip

\subsection*{Acknowledgment} The first author is supported by
a Summer Undergraduate Research Fellowship at Caltech. The
second author is supported by NSF grants
DMS-0901221, DMS-1007207, DMS-1201512, PHY-1205440.

\end{document}